\input amssym.def 
\magnification = \magstep2
\baselineskip=15 true pt
\rightskip=-0.1cm

\def\lf{\ \hfil\break}       
\def\LF{\medskip\noindent}   
\def\cl{\centerline}
\def\id{{\rm id}\,}
\font\cmrX=cmbx10 scaled \magstephalf
\def\BR{\Bbb R}
\def\BH{\Bbb H}
\def\BS{\Bbb S}
\font\eu=eusb10
\def\SC{\hbox{\eu C}}

\cl{\cmrX Riemannian Center of Mass and so called karcher mean }
\lf
\cl{Hermann Karcher, Bonn}
\LF
{\narrower\noindent {\bf Abstract.} The {\it Riemannian center of mass} was constructed
in [GrKa] (1973). In [GKR1, GKR2, Gr, Ka, BuKa] (1974-1981) it was successfully 
applied with more refined estimates. Probably in 1990 someone renamed it without 
justification into {\it karcher mean} and references to the older papers were
omitted by those using the new name. As a consequence newcomers started to reprove 
results from the above papers. -- Here I explain the older history.
\par}

\LF
The Euclidean center of mass is an {\bf affine} notion. I will use discrete mass points
rather than mass distributions, for simpler wording. Let $ p_i \in \BR^n$ be points
of weight $m_i$ with $\sum m_i =1$ (for convenience). The center is defined, without 
reference to a metric, as:
 \vskip5pt
\cl{$\SC(p_i) :=\sum m_i\cdot p_i .$} \lf
This formulation does not carry over to Riemannian manifolds. A small change improves
the situation. We define a vector field on $ \BR^n$ by averaging the vectors from
an arbitrary point $x\in \BR^n$ to the mass points $p_i$:
 \vskip5pt
\cl{$V(x) :=\sum m_i\cdot (p_i- x). $} \lf
This vector field points from every point $x$ to the center, $\SC(p_i) = x + V(x)$,
and it generalizes to a Riemannian manifold $M$:
 \vskip5pt
\cl{$V(x) :=\sum m_i\cdot \exp^{-1}_x p_i. $} \lf
Under the assumptions of [GrKa] we proved that the total index of the vector field 
$-V$ on a convex ball containing the $p_i$ is $1$ and the covariant differential $-DV$ is
close to $\id$ so that the index at each zero is also $1$. Therefore $V$ has a
unique zero in the convex ball. We called it the {\it Riemannian center of mass}.
The estimate on the covariant differential is also used to prove that the
geodesic $ \gamma:[0,1]\to M,\ \gamma(0)=x,\ \gamma'(0)=V(x)$ ends closer to the center 
than $x$. This geodesic can be thought of as one Euler step for computing the
integral curve of $V$ starting at $x$. All the papers quoted above explicitly say
how much closer, at least, one such Euler step gets to the center.
\lf
In the Euclidean case the vector field $-V$ is the gradient of the convex function
$f(x):=0.5\sum m_i\cdot |p_i- x|^2$ for {\bf any} choice of the Euclidean metric
$|.|^2$ on $\BR^n$. In the Riemannian case our estimate on $DV$ says that the corresponding
function is convex [Ka], but it is not surprising that {\bf no} estimate in later papers 
uses the {\bf values} of that function.
\lf
The hyperbolic space $\BH^n$ and a hemisphere of $\BS^n$ are not linear spaces, but they 
have standard isometric embeddings into $\BR^{n+1}$. One can {\bf explicitly} compute the
center in $\BR^{n+1}$ and project it back into $\BH^{n}$ or $\BS^{n}$ to get an isometry 
invariant explicit center in these Riemannian manifolds. I do not know how old this 
construction is. It was known in the seventies, but it gives no clue what to do if the 
sectional curvature is not constant.
\lf
We found the center of mass a surprisingly effective tool. In [GrKa] we
define the diffeomorphism that conjugates two $ C^1$-close group actions with
a single center of mass application. In [GKR1] we improve an almost homomorphism
of compact groups by an iteration where each step is one center of mass application.
In [GKR2] we improve previous results to get dimension independent estimates by using 
a Finsler metric rather than the Riemannian metric on the orthogonal group. This
improvement emphasizes the affine background of the center because the biinvariant
connection on the group is a metric connection for the Riemannian and for the Finsler metric.
The exponential map depends only on the connection, the center is therefore
the same for both metrics, but the needed convex sets are much larger and the estimates 
much better in the Finsler case. [Ka] was written for my hosts at the Courant
Institute. Therefore all Jacobi field estimates are proved via differential inequalities,
i.e. with tools from analysis rather than from differential geometry. The main result is in
the title. [BuKa] explains in the first five chapters Gromov's proof and in the remaining
three supplies the needed differential geometric tools, written for a broader audience.
Chapter 6 contains Jacobi field estimates in the Riemannian case and how they control
geodesic constructions. Chapter 7 repeats this for Lie groups, demonstrates the advantage
of non-Euclidean norms, converts a power series solution of the Jacobi equation into
a differential equation for the Campbell-Hausdorff formula. Chapter 8 contains a detailed 
exposition of the center of mass: it begins with the Euclidean case, covers the
Riemannian case, compact Lie groups, finally nilpotent Lie groups (no size restriction
on the support of the masses) and describes applications.
\lf
Minima of convex functions have been used much earlier by Cartan in mani\-folds of negative
curvature. Our construction initially aimed at positively curved manifolds and could not have been
done before Rauch introduced Jacobi field estimates as a tool in Riemannian geometry.
I mention convex functions in [Ka] since our estimates imply convexity results. But I say
on p.511 that the vector field is more important and that $\sum m_i\cdot(1-\cos(d(x,p_i))$
gives, on $\BS^n$, a better curvature adaption then 
$0.5\sum m_i\cdot d(x,p_i)^2$. Such modified distance functions are a main tool in [ChKa].
\LF
Then, around 2000, a renaming caught on. The earliest use of {\it karcher mean} I know of
is in a 1990 paper by W.S.Kendall, also cited in [Af] as the earliest case.
I know of no other such renaming decision. In the {\it karcher mean} literature the older papers
are essentially ignored. Many of our estimates are reproved
and the advantage of Finsler metrics on the orthogonal groups has not yet been noticed.
I think it is fair to say that a substantial amount of damage was caused by the renaming.
\lf
Wikipedia talks about a generalization of centroids to metric spaces {\bf if} the minimizer
of the function $f(x):= \sum m_i\cdot d(x,p_i)^2$ {\bf exists and is unique} and continues: 
\lf {\it
``It is named after Maurice Fr\'echet. Karcher means are a closely related construction
named after Hermann Karcher.'' }
\lf
Who believes that Fr\'echet would have been proud of this naming?
\vskip6mm

\centerline{{\bf Bibliography} (chronological)}\LF
On the center of mass: \LF
[GrKa] Grove,K., Karcher,H.: How to conjugate $ C^1$-close group actions. 
\hglue1.25cm Math.Z. {\bf 132}, 11-20 (1973).
\lf
[GKR1] Grove,K., Karcher,H.,Ruh,E.: Group actions and curvature. Inven-
\hglue1.25cm tiones math. 
{\bf 23}, 31-48 (1974).
\lf
[GKR2] Grove,K., Karcher,H.,Ruh,E.: Jacobi fields and Finsler metrics on 
\hglue1.25cm  compact Lie groups with an application to the differentiable pin-
\hglue1.25cm  ching problem. Math. Ann. 
{\bf 211}, 7-21 (1974).
\lf
[Gr]\hskip18ptGrove,K.: Center of Mass and G-local triviality of G-bundles.
Proc. \hglue1.25cm AMS {\bf 54}, 352-354 (1976).
\lf
[Ka]\hskip18ptKarcher,H.: Riemannian center of mass and mollifier smoothing.
\lf \hglue1.25cm CPAM {\bf 30}, 509-541 (1977).
\lf
[BuKa] Buser,P., Karcher,H.: Gromov's almost flat manifolds. Ast\'erisque 
 \hglue1.25cm {\bf 81}, 1-148 (1981).
\LF 
Using curvature adapted distance functions:
\LF
[JoKa] Jost,J., Karcher,H.: Geometrische Methoden zur Gewinnung von 
a-\hglue1.2cmpriori-Schranken f\"ur harmonische Abbildungen. Manuscripta math. 
\hglue1.25cm {\bf 40}, 27-77 (1982).
\lf
[ChKa] Karcher,H.:  Riemannian comparison constructions.
In Chern,S.S. \lf \hglue1.25cm (editor): Global Differential Geometry.
MAA Studies in Mathemat-
\hglue1.25cm ics {\bf 27}, 170-222 (1989).
\LF 
With four pages of history:
\LF
[Af]\hskip18ptAfsari,B.: Riemannian $L^p$ Center of Mass: Existence, 
Uniqueness \hglue1.26cm and Convexity. Proc. AMS {\bf 139}, 655-673 (Feb 2011).
\LF

Hermann Karcher \hfill June 2014

Mathematisches Institut der Universit\"at

D-53115 Bonn

unm416@uni-bonn.de
\bye